\documentstyle[11pt,amsfonts]{article}
\def\cancel#1#2{\ooalign{$\hfil#1\mkern1mu/\hfil$\crcr$#1#2$}}
\newcommand{\ba}{\begin{eqnarray}}
\newcommand{\na}{\end{eqnarray}}
\newcommand{\ban}{\begin{eqnarray*}}
\newcommand{\nan}{\end{eqnarray*}}
\newtheorem{The}{Theorem}[section]
\newtheorem{Cor}[The]{Corollary}
\newtheorem{Pro}[The]{Proposition}
\newtheorem{Ex}[The]{Example}
\newtheorem{R}[The]{Remark}
\newtheorem{Le}[The]{Lemma}
\newtheorem{D}[The]{Definition}

\def\endproof{\relax\ifmmode\expandafter\endproofmath\else
  \unskip\nobreak\hfil\penalty50\hskip.75em\hbox{}\nobreak\hfil\bull
  {\parfillskip=0pt \finalhyphendemerits=0 \bigbreak}\fi}
\def\endproofmath$${\eqno\bull$$\bigbreak}
\def\bull{\vbox{\hrule\hbox{\vrule\kern3pt\vbox{\kern6pt}\kern3pt\vrule}\hrule}}
\def\dirac{\mathpalette\cancel D}
\def\dir{\mathpalette\cancel\partial}

\begin{document}
\title{Seiberg-Witten-Floer Theory for Homology 3-Spheres}
\author{Bai-Ling Wang}
\date{Feb. 10, 1996}
\maketitle

\begin{abstract}
We give the definition of  the Seiberg-Witten-Floer homology group for a
homology 3-sphere.  Its Euler characteristic number is a Casson-type 
invariant. 
For a four-manifold with boundary a homology sphere, a relative 
Seiberg-Witten invariant is defined taking values in the Seiberg-Witten-Floer 
homology group, these relative Seiberg-Witten invariants are applied
to certain homology spheres bounding Stein surfaces.

\end{abstract}

\section{Introduction}

Since Seiberg and Witten introduced a new gauge theory through the monopole 
equations \cite{W2}, the power of these new invariants arising from 
the monopoles and the subtlety they capture through the scalar curvature 
has been used to prove  the ``Thom" conjecture and 
progress has been made 
on the $``\frac {11}{8}"$ conjecture \cite{KM} \cite{Fru}.
Furthermore, as conjectured by Witten and supported by physicists'
calculations, the Seiberg-Witten invariants appear to be
  equivalent to Donaldson's
polynomial invariants. 

 In recent years progress on instanton theory has been rapid
through the study of  the relative Donaldson invariant for 4-manifolds with
cylindrical ends, where Floer instanton homology enters very naturally. 
It is reasonable  to expect that Seiberg-Witten theory for 3-manifolds could
play a similar role. In this 
paper, we give a mathmatically rigorous definition of an analogue of the Floer
homology group for a homology 3-sphere where in this case any
line bundle has trivial first Chern class.  
 For a 3-manifold with nontrivial first Betti number,
Marcolli \cite{Mar2} constructs the Seiberg-Witten-Floer theory,
which is independent of metric and perturbation, therefore,
a topological invariant. In our case, since the reducible solution,
 whose spinor is vanishing, 
cannot be perturbed away,  more care must be taken for the reducible
solutions. For those metrics whose ordinary Dirac operator have trivial 
kernel, then the reducible solution as a critical point 
is isolated and unique, we can define the Seiberg-Witten-Floer homology
by removing this reducible solution. 
In this present paper, we analyse the asymptotic behaviour of the
gradient flows connecting two critical pionts. Using the
gluing arguments, we show that the  Seiberg-Witten-Floer homology
(removing the reducible solution) is well-defined.
 The subject of metric (perturbation) dependence
is  interesting. In concurrent joint work with M. Marcolli, 
we study  equivariant
 Seiberg-Witten-Floer theory which is metric independent up to 
{\it index-shifting}. Interestingly, not only are there  
Seiberg-Witten invariants  for two critical points with relative
index 1, but also  Seiberg-Witten invariants
 for two critical points with relative index 2, which essentially   
measure the interactions of the irreducible critical points
with the unique reducible one (see \cite{MW} for details). 

Our main strategy is to study the Seiberg-Witten equations near the critical
points using hyperbolic, non-linear equations. We give
a detailed analysis to derive exponential decay estimates for the
solutions and use these decay estimates to get a useful gluing theorem
analogous to  
 {\em Taubes' constructions} in instanton theory. These gluing 
theorems reflect the fact that  Floer's work for instanton homology can
be adapted to define the Seiberg-Witten-Floer homology invariant. 
It is well-known now that the Floer homology group has a non-trivial 
cup-product (the so-called quantum cup-product), 
which has encoded the rich
information from quantum string theory. 
One may expect these also from the Seiberg-Witten-Floer homology group.
(In preparation is the case of 
a circle bundle  over a Riemann 
surface.)
This is likely to
be a very interesting topic for the future in the field of topological 
quantum field theory and topological quantum gravity. The mathematical
part of this is Taubes' astonishing result in \cite{Tau1}
 which asserts that the 
Seiberg-Witten invariant for a symplectic manifold $X$  equals
 the Gromov invariant for a homology class in $H_2 (X, {\Bbb Z})$ associated with
the $Spin^c$ structure for the Seiberg-Witten invariant.

\section{Seiberg-Witten equations on homology 3-sphere $Y$ and its cylinder
$Y\times {\Bbb R}$}
\subsection{Seiberg-Witten equations on $Y \times {\Bbb R}$}

In this subsection, 
we review the Seiberg-Witten equations on a compact, connected,
closed, oriented, homology 3-sphere $Y$ and its cylinder $Y \times {\Bbb R} $, with a
 $Spin^c$ structure on $Y \times {\Bbb R} $ given by pulling back 
a $Spin^c$ structure on $Y$. 
The two $Spin^c$ bundle $S^+$ and $S^-$ can be identified via the Clifford
multiplication by $dt$,  with the 
$Spin^c$ bundle $S$ on $Y$. Since $Y$ is a spin manifold, we can write 
the determinant bundle of $S$ as $L^2$ for some line bundle $L$ on $Y$.

The Seiberg-Witten equations for $Y \times {\Bbb R }$  are the equations for
a unitary $U(1)$-connection $A$ on $L$ and a spinor section
 $\psi\in \Gamma (S)$:
\ba
&&F^+_A = \frac{1}{4}<e_ie_j\psi, \psi>e^i\wedge e^j\nonumber\\[1mm]
&&\dirac_A (\psi ) =0 
\label{SW4}
\na 
where $F^+_A$ is the self-dual part of the curvature of $A$, $\{e_i\}_{i=1}^4$
 is the 
 orthonormal frame for $T(Y\times {\Bbb R})$ 
acting on spinors by  Clifford multiplication. We adopt the usual convention
for the Clifford algebra:
$e_ie_j +e_je_i = -2 \delta_{ij}$, with $\{e^i\}_{i=1}^4$
its dual basis elements for $T^*(Y\times {\Bbb R})$. We
denote by $<,>$ the hermitian 
inner product and sometimes write $\overline \psi _1 \psi _2 
= <\psi_1, \psi_2>$. Here    
$\dirac_A $ is the Dirac operator on $Y\times {\Bbb R}$
 associated with the connection $A$. 

We say that a pair $(A, \psi)$ is in temporal gauge if the $dt$-component of
$A$ vanishes identically. In this temporal gauge, $(A, \psi )$
on $Y \times {\Bbb R}$ can be written as a path $(A(t), \psi(t) )$ on $Y$ obtained
by restricting $(A, \psi)$ to the slices $Y \times \{t\}$.  Then the 
Seiberg-Witten equations (\ref{SW4}) read as follows (\cite{KM}, \cite{Marcolli}, \cite{Wang2}). 

\ba
&&\frac{\partial A }{\partial t} = *F_A - q (\psi)\nonumber\\[1mm]
&&\frac{\partial \psi }{\partial t} =  \dir _A (\psi)
\label{SW4t}
\na
where $ q( \psi ) $ is the quadratic function of $\psi$, in local coordinates,
$$
q(\psi) = \frac{1}{2}<e_i.\psi, \psi> e^i,
$$
and $*$ is the Hodge star operator on $Y$, $\dir_A$ is the Dirac operator on $Y$
twisted with a time-depedent connection $A$.

As we work in temporal gauge, the permitted gauge transformations are
constant with respect to the $t$-direction and we denote this gauge group by
${\cal G} = Map(Y, U(1))$. It is easy to see that
the  equations (\ref{SW4t}) are invariant under
 the gauge group $\cal G $
$$(A, \psi) \mapsto (A - id\tau, e^{i\tau}\psi).$$  
 
As first noted in \cite{KM}, the Seiberg-Witten equations (\ref{SW4t}) 
are a gradient flow equation on $\cal A$ (the space of pairs $(A, \psi)$) of
a functional $C$. Before we introduce this functional, let us give an
 appropriate
Sobolev norm on the configuration space $\cal A$. We choose to work 
with $L^2_1$-connections $A$ and $L^2_1$-sections of the associated $Spin^c$
bundle, and the gauge transformations in $\cal G$ consist of
$L^2_2$-maps from $Y$ to $U(1)$, these gauge transformations are at least
continuous by the Sobolev embedding theorem. 
With this metric on $\cal A$, the tangent space of $\cal A$ is
the space of $L_1^2$-sections of $\Omega (Y, i {\Bbb R}) \oplus \Gamma (S)$.  The 
inner product on this tangent space is the $L^2$-product on the one forms and 
twice the real part of the $L^2$-hermitian product, that is, for the tangent
vectors $a, b \in \Omega (Y, i {\Bbb R})$, the $L^2$ inner product is
given by 
$$
\langle a, b\rangle_{L^2} = -\int_Y a\wedge *b
$$
where $*$ is the complex-linear
 Hodge-star operator on $Y$, for the tangent vectors
$\psi_1, \psi_2 \in \Gamma (S)$, the $L^2$ inner product is given by
$$
\langle \psi_1, \psi_2 \rangle_{L^2} = \int _Y( \overline \psi_1 \psi_2 +   
 \overline \psi_2 \psi_1 )
$$
We drop the subscript $L^2$ when there  is no danger of
confusion in the following
definition of the functional.

\begin{D}
Fix a $U(1)$-connection $B$ on $L$, define the functional $C$ as
\ba 
C(A, \psi) = \frac{1}{2} \int_Y (A-B) \wedge F_A +  \langle \psi, \dir _A \psi
\rangle d vol(Y)
\label{func}
\na
\end{D}

For our ${\Bbb Z}$-homology sphere, this functional can be reduced to the 
quotient space ${\cal B} = {\cal A} /{\cal G}$
 (which is a Hausdorff space by the
following lemma). Denote by ${\cal B}^* $ the space of the irreducible pairs
$(A, \psi)$ (where $\psi \ne 0$) modulo the action of the gauge group $\cal G$.
The following lemma is a standard result of elliptic regularity.

\begin{Le}
(1) The  gauge group $\cal G$ acts smoothly on $\cal B$. The quotient 
${\cal B}^*$
is a smooth Hilbert manifold with tangent spaces
$$
T_{[A, \psi]}({\cal B}^* ) = \{ (\alpha, \phi )\in \Omega (Y, i {\Bbb R})\oplus 
\Gamma (S) | d^* \psi = \overline \psi \phi -\overline \phi \psi \}. 
$$\\
(2) The charts for ${\cal B}^* $ near $[A, \psi]$ are given by 
\[
T_{[A, \psi]}({\cal B}^* ) \longrightarrow {\cal B}^* ; \qquad \quad (a, \phi)\mapsto
[A +a, \psi +\phi ]
\]
\label{TB}
\end{Le}

\begin{Pro}(\cite{KM}\cite{MST})
Fix an  interval $[0, 1]$, if $[A, \psi]$ is the solution of the 
Seiberg-Witten equations (\ref{SW4t}), then as a path in $\cal B $,
$[A, \psi] $ satisfies the gradient flow equation 
$$
\frac {\partial (A, \psi) }{\partial t} = \nabla C _{[A, \psi]}
$$
where $\nabla C _{[A, \psi]} = (*F_A - q(\psi), \dir_A (\psi))$
 is the $L^2$-gradient vector field on
$\cal B$ defined by the functional $C$.
\end{Pro}

\subsection{Seiberg-Witten equations on $Y$ and critical points for $C$}

We have shown that the critical points of $C$ are solutions of 
the Seiberg-Witten 
equations on $Y$,
\ba
*F_A& =& q(\psi)\nonumber \\[2mm]
\dir_A(\psi)& =&0
\label{SW3}
\na

For the Seiberg-Witten equations (\ref{SW3}) on a homology 3-sphere $Y$, we
find an interesting phenomenon, that is the reducible solutions cannot
be perturbed away by changing the  metric or perturbing the
 curvature equations.
However, for the metric whose oridinary Dirac operator has
trivial kernel, the reducible solution as a critical point in ${\cal B}$
is non-degenerate, hence isolated. Note that for the homology sphere,
 there is only one reducible critical point.  

\begin{Le}
For the metric $g$ on $Y$, whose ordinary Dirac operator $\dir$ has
trivial kernel, then the only reducible critical point  $[0, 0]$,
  the trivial solution for (\ref{SW3}), is isolated and nondegenerate.
\label{iso}
\end{Le}

\begin{proof}
Since our $Y$ is a $\Bbb Z$-homology 3-sphere, the reducible solution for
(\ref{SW3}) is $(A, 0)$ where $F_A =0$, that means that 
$A$ is abelian flat connection, 
determined by the $U(1)$-representation for the fundamental group $\pi_1(Y)$. 
The connection $A$ must be in the
 orbit through  the trivial connection on $L$. 
The isolated property can be proved using the  Kuranishi model for the
 reducible
solution $[0, 0]$, here we use the perturbation theory by expanding the solution
near $[0, 0]$. First, note that if a solution $[A, \psi]$
 of (\ref{SW3}) is sufficiently
close to $[0, 0]$, then $ [A, \psi]$ obeys the following equations
\ban
&& d^* A =0 \\[1mm]
&& *dA = \overline \psi \sigma \psi\\[1mm]
&& \dir (\psi) + A. \psi =0
\nan 
where $\sigma = c(e_i) e^i$, $c(e_i)$ is the Clifford multiplication of 
$\{e_i \}_{i=1}^3$. We can write $(A, \psi)$ as
\ban 
&&A= \epsilon a_1 + \epsilon ^2 a_2 + \epsilon^3 a_3+\cdots\\[1mm]
&&\psi = \epsilon \psi _1 + \epsilon^2 \psi_2 +\epsilon^3 \psi_3+\cdots\\[1mm]
\nan
Then from the Dirac equation and $ ker \dir =0$, we know that $ \psi _1 =0$,
use this fact and the curvature equation, we know that $a_i$ $ (i=1, 2, 3)$ 
are colsed and co-closed, for a homology 3-sphere, this leads to 
$a_i (i=1, 2, 3)$ must be zero. Repeat this procedure, we get $[A, \psi]$
is zero. 
The non-degenerate property for $[0, 0]$ is equivalent to the isolated property
in some sense. The Hessian operator at $[0, 0]$ is the map,
$$
K: \qquad T_{[0, 0]}({\cal B}) \longrightarrow T_{[0, 0]}({\cal B}) $$
$$ (a, \phi ) \mapsto (*da, \dir (\phi) )$$
where $d^* a =0$, it is easy to see that $ker K =0$.
\end{proof} 
  
Let ${\cal M} _Y$ denote the set of the critical points of $C$ on $\cal B$, 
i.e. the set of solutions of (\ref{SW3}) modulo the gauge group $\cal G$
action.

For $[A_0, \psi_0] \in {\cal M} _Y, \psi_0 \ne 0$, the Hessian  operator $K$
is the following operator, it is the linearization of (\ref{SW3}) on 
${\cal B}^*$,
$$
K_{[A_0, \psi_0]}: \qquad T_{[A_0, \psi_0]}({\cal B}^*) \longrightarrow 
T_{[A_0, \psi_0]}({\cal B}^*) $$
\[  \left(\begin{array}{c}a \\ \phi\end{array} \right) \mapsto 
\left( \begin{array} {cc} *d & -Dq_{ \psi_0}\\  .\psi_0 &  \dir _{A_0}
\end{array}\right)\left( \begin{array}{c}a\\ \phi \end{array}\right) \]
where $Dq_{\psi_0}$ is the linearization of $ q$ at $\psi_0$,
and $\psi_0$ acts on $a$ by $a. \psi_0$.

In general, $K$ may not be non-degenerate at $[A_0, \phi_0]$, but it has index
zero, this can be seen from the fact 
 that $K$ is a compact perturbation of the 
Dirac operator and 
\[ \left( \begin{array}{cc} 0 & d^* \\
                            d & *d\end{array}\right): \quad
         \Omega^0(Y, i{\Bbb R}) \oplus \Omega^1(Y, i{\Bbb R}) \rightarrow  
        \Omega^0(Y, i{\Bbb R}) \oplus \Omega^1(Y, i{\Bbb R})  
\]
Both have index zero on three manifolds. This means that after a generic
perturbation, $K$ would have non-degenerate critical points set ${\cal} M_Y$.
Before introducing the perturbation, we study the geometry of the moduli
space ${\cal M}_Y$ from the deformation complex at $[A_0, \phi_0]$:

\ba
0\to \Omega^0(Y, i{\Bbb R}) \stackrel{G}\to \Omega^1(Y, i{\Bbb R}) \oplus \Gamma (S)
 \stackrel{L}\to   \Omega^1(Y, i{\Bbb R}) \oplus\Gamma (S)/(G(\Omega^0(Y, i
{\Bbb R}))
\label{deform}
\na
where $ G$ is the infinitesimal action of the gauge group $\cal G$,
$$ \tau \mapsto (-d\tau, \tau \psi)$$
and the map $ L$ is the linearization of (\ref{SW3}) on $\cal C$ at 
$(A_0, \psi_0)$,
\[ 
\left( \begin{array} {c} a \\ \phi \end{array}\right)\mapsto
\left( \begin{array} {cc} *d & -Dq_{ \psi_0}\\  .\psi_0 &  \dir _{A_0}
\end{array}\right)\left( \begin{array}{c}a\\ \phi \end{array}\right) \] 

The Zariski tangent space of ${\cal M}_Y$ at $[A_0, \phi_0]$ is given by 
\[
ker(L) /Im (G).\]

In order to see that $ {\cal M}_Y^* = {\cal M}_Y \setminus \{[0, 0]\}$ 
is a 0-dimensional
smooth submanifold in ${\cal B}^*$, 
we need to consider
a generic perturbation for (\ref{SW3})  by proving
$ker (L)= Im (G)$. 
In instanton theory, there is a generic metric theorem. It states that
for a generic metric (that is lying in an
open dense subset in the space of metrics) the
instanton moduli space is a smooth manifold. There is no such 
generic metric theorem for us. But as in \cite{KM} \cite{W2}, we can perturb
the curvature equation by adding a suitable 1-form.

Now consider the perturbed functional $C'$ as follows,
\[
C'(A, \psi) = \frac {1}{2}\int_Y (A-B) \wedge( F_A -2 d\alpha ) +  \langle \psi, \dir _A \psi
\rangle d vol(Y)
\label{func'}
\]
where $\alpha$ is a purely imaginary 1-form on Y, as we remark below, this
is the only possible perturbation for the curvature equation.

Then the perturbed Seiberg-Witten equations on $Y$ (the solutions are
the  critical points for $C'$) are given by
\ba
&& *F_A =  q(\psi ) + *d \alpha \nonumber\\[2mm]
&& \dir_A (\psi) =0
\label{SW3'}
\na
while the gradient flow equation (the perturbed Seiberg-Witten equations on 
the cylinder $ Y \times {\Bbb R}$) is given by
\ba
&&\frac{\partial A }{\partial t} = *F_A - q (\psi) -*d\alpha \nonumber\\[1mm]
&&\frac{\partial \psi }{\partial t} =  \dir _A (\psi)
\label{SW4t'}
\na
\begin{R}
We point out that $ q(\psi ) $ is co-closed for $\psi$ satisfying the
Dirac equation $\dir_A (\psi) =0$, therefore the perturbing form in (\ref{SW3'})
(\ref{SW4t'}) must be also co-closed by the Bianchi identity $d^* (*F_A)=0$.
For a homology 3-sphere, it is $ *d\alpha$ where $\alpha$ is any 
imaginary-valued one form on $Y$.
\end{R}

We denote the moduli space  ${\cal M}_{Y, \alpha}$ for the solution space of 
(\ref{SW3'}) modulo $\cal G$. 
\begin{Pro}
For any metric $g$ on $Y$, there is an open dense set in 
$*d(\Omega_{L^2_2}^1(Y, i{\Bbb R}))$, such that ${\cal M}_{Y, \alpha}^*$ is 
a smooth 0-dimensional manifold and the only reducible solution $[\alpha, 0]$
is isolated for  any $*d\alpha$ in this set.
\label{generic}
\end{Pro}
\begin{proof}
This is an application of the 
 Sard-Smale theorem as \cite{KM} \cite{Mar2}. First we 
construct the parametrized moduli space,
\ban \begin{array}{lcl}
{\cal M} &=& \bigcup _{\delta =*d\alpha} {\cal M}_{Y, \alpha}\\[2mm]
      &\subset & ({\cal A} \times *d(\Omega_{L^2_2}^1(Y, i{\Bbb R})))/{\cal G}
\end{array}
\nan
We write ${\cal M}^*, {\cal M}_{Y, \alpha}^*$ for the moduli spaces of irreducible 
solutions.
Let $(A_0, \psi_0, \delta)$ be the solution of (\ref{SW3'}) with $\psi \ne 0$.
The linearization of the equations at this point is the following operator:
\[
P: \qquad \Omega^1(Y, i{\Bbb R}) \oplus\Gamma (S) \oplus *d(\Omega^1(Y, 
i{\Bbb R})) 
\rightarrow \Omega^1(Y, i{\Bbb R}) \oplus\Gamma (S)/(G(\Omega^0(Y, i{\Bbb R})) \]
\[
\left(\begin{array} {c} a \\ \phi \\ \eta \end{array}\right)\mapsto
\left( \begin{array} {ccc}  *d & -Dq_{\psi_0} & -Id \\
       .\psi_0 &  \dir _{A_0} & 0 \end{array}\right)
\left( \begin{array}{c}a\\ \phi \\ \eta \end{array}\right) \]
To estabilish surjectivity for $P$, it is sufficient to show that no element
$(b,  \rho )$ in the range can be $L^2$-orthogonal to the image of $P$.
If such $(b,  \rho )$ is orthogonal, then by varying $\eta$ alone,
one can get  that $b$ is 0. Varying $\phi$ alone, $\rho$ must be in the kernel
of $\dir_{A_0}$, then varying $a$ alone, one can see $\rho$ is also 0,
otherwise since
being a solution of the Dirac equation $\dir_{A_0} =0$, $\psi$ and $\rho$
cannot vanish on an open set, note that, modulo the gauge action in 
the range, we can construct an $a$ such that $(0, \rho)$ is not orthogonal
to $P(a,0,0)$.

Hence by the implicit function theorem, we know that ${\cal M}^*$
is a smooth manifold, now apply the  Sard-Smale theorem to the projection 
map: 
\[ \pi : \quad {\cal M}^* \longrightarrow *d(\Omega^1(Y, i{\Bbb R}) \]
we get there is a generic $\delta = *d\alpha$ such that  
${\cal M}_{Y \alpha}^*$ is a smooth manifold, by the index calculation, it
is 0-dimensional. 

Consider the isolated property of the only reducible point $[\alpha ,0]$.
By an argument  
similar to the proof of Lemma \ref{iso},  the perturbing form $*d\alpha $ 
should satisfy $ker \dir_{\alpha}=0$, this is a co-dimension one condition
for $\alpha$. The existence of an open dense set of such 
$\alpha $ still holds.  
\end{proof}

\section{Properties of various moduli spaces}

\subsection{The compactness of the moduli space}

There is a great simplicity in Seiberg-Witten gauge theory, 
that is, the moduli space of Seiberg-Witten
equations is always compact. In this subsection,
we will show that ${\cal M}_{Y, \alpha}$ is compact and, 
since the reducible point in ${\cal M}_{Y, \alpha}$ is isolated,
for a generic sufficiently small 
perturbation, ${\cal M}_{Y, \alpha}^*$ is a set with finitely many points,
 consisting
of the irreducible, non-degenerate critical points of $C'$. Moreover, we can 
also prove the moduli space of the Seiberg-Witten equations (\ref{SW4t'})
connecting the two critical points in ${\cal M}_{Y, \alpha}$ is also compact.

\begin{Pro}
${\cal M}_{Y, \alpha}$ is sequentially compact, that is, for any sequence
$\{(A_i, \psi_i)\} $ of the solution (\ref{SW3}) or (\ref{SW3'}), there
is a subsequence (up to $L_2^2$ gauge tranformations) converging to a
solution in $C^\infty$-topology. Therefore, ${\cal M}_{Y, \alpha}$ contains
only finite points in $ \cal B$. 
\end{Pro} 
This follows from a priori bounds for 
any solution of ($\ref{SW3}$). 
 There are several proofs on the
compactness of the Seiberg-Witten moduli space in the literature such as 
\cite{KM} \cite{T2} \cite{Marcolli}. 
We omit the proof of this proposition.

Now we are in a position to discuss the moduli space of the connecting orbits.
The functional $C$ or $C'$ has a nice property, since it satisfies the 
Palais-Smale condition. 

\begin{Le}
For any $\epsilon > 0$, there  is $\lambda > 0$ such that if $[A,\psi]
\in {\cal B} $ has the $L^2_1$-distance at least $\epsilon $ from all
the critical points in ${\cal M}_{Y, \alpha}$, then 
\[
\| \nabla C'_{[A, \psi]}\|_{L^2} > \lambda 
\]
\label{ps}
\end{Le}
A
similar result was obtained in \cite{MST} for three manifolds $S^1 \times 
\Sigma $, where $ \Sigma $ is a Riemann surface of genus $g>1$
 whose determinant line bundle
for the $Spin^c$ bundle is pulled back from a line bundle on the Riemann
surface which has degree $2-2g$.

\begin{proof}
Suppose there is a sequence $(A_i, \psi_i)$ in $\cal B$ 
whose $L_1^2$-distance are at least $\epsilon $ from all 
the critical points in ${\cal M}_{Y, \alpha}$ for which
\[
\| \nabla C'_{[A_i, \psi_i]}\|_{L^2} \to 0 \qquad \qquad \hbox{as}\ i\to \infty 
\]
Then as $i\to \infty$, we have
\[ \begin{array}{r}
\| * F_{A_i} - q(\psi_i) - * d \alpha \|_{L^2} \to 0 \\
\| \dir_{A_i} (\psi _i) \|_{L^2} \to 0 
\end{array}
\] 
where $*d\alpha$ is sufficiently small. 
This means that there is a contant $C>0$ such that
\[
\int_Y |* F_{A_i} - q(\psi_i) - * d \alpha |^2 + |\dir_{A_i} (\psi _i)|^2
< C
\]
Resorting to the  Weitzenbock formula for $\dir_{A_i}^2$ and for 
$*d\alpha$ sufficiently
 samll, the above inequality reads as 
\[
\int_Y |F_{A_i}|^2+|q(\psi_i)|^2 + \frac s2 |\psi _i|^2 + 2 |\nabla_{A_i}\psi_i|
^2  < 2 C 
\]
It follows that $ \| \psi _i\|_{L^4}, \|F_{A_i}\|_{L^2},  
\|\nabla_{A_i}\psi_i \|_{L^2}$ are all bounded independent of $i$.
Now use the standard elliptic argument as in the proof of the compactness, 
we know that there is a subsequence converging in $L^2_1$-topology to a
solution of (\ref{SW3'}), this contradicts with the assumption 
$ (A_i, \psi_i)$ in $\cal B$ 
whose $L_1^2$-distance are at least $\epsilon $ from all 
the critical points in ${\cal M}_{Y, \alpha}$.
\end{proof}

\begin{D}
A finite energy solution of  (\ref{SW4t'}) on $ Y\times {\Bbb R}$
 is the solution $ [A(t), \psi(t)]$  whose square of the energy 
\[
\int^{+\infty}_{-\infty} \|\nabla C'([A(t), \psi(t)])\|^2_{L^2(Y)}
\]
is finite. 
\end{D}

Choose a sufficiently small $\epsilon$ for all the finite points 
in ${\cal M}_{Y, \alpha}$ such that Lemma \ref{ps} holds, we label the 
nondegenerate, irreducible  critical 
points of $C'$ as $ x_1, \cdots , x_N$.  
Denote by $U(x_i, \epsilon)$ the $L_1^2$-open ball of radius
 $\epsilon$. Let $\gamma (t) = [A(t), \psi (t)]$
is the gradient flow for $C'$ with finite energy. From Lemma \ref{ps}, 
one can see
\[
\{ t\in {\Bbb R }| [A(t), \psi (t)] \in U(x_i, \epsilon)\ \hbox {for any}\ i \}
\]
has total length finite and has only finite intervals running between 
different  $U(x_i, \epsilon)$.
 
Therefore, for sufficiently large $T>>1$, $\gamma |_{[T, \infty )} \subset 
 U(x_i, \epsilon)$ for some $x_i$. 

Similarly,  $\gamma |_{[-T, -\infty )} \subset  
U(x_j, \epsilon)$ for some $x_j$.
Actually, one can use the finite energy to get the limits of $\gamma (t) $
which are $x_i$, $x_j$ respectively. In the next subsection, we will
prove that $\gamma (t) $  exponentially decays to the limits.

We denote the moduli space of finite energy
 solutions for (\ref{SW4t'}) which connect
the two critical points $a, b$ by 
\[ 
{\cal M} (a, b) = 
\left\{
x(t)\left| \begin{array}{l}
(1)\ \hbox{$x(t)$ is the finite energy solution of} \\[2mm]
    \quad \hbox{ (\ref{SW4t'}) modulo the gauge group $\cal G$,}
\\[2mm]
(2)\ \lim_{t \to +\infty}x(t)= a ,\\[2mm]
(3)\ \lim_{t \to -\infty}x(t)=b. 
\end{array} \right. \right\}
\]

\begin{R}
Since any Seiberg-Witten monopole on $S^1 \times Y$ is invariant under
the rotation action of $S^1$, there is no closed  non-constant gradient flow
connecting the same critical point, that is,
${\cal M}(a, a) $ is just one piont  $\{a\}$.
One can also  see this from the following
useful length equality:
\[
\begin{array} {ll}
& C' ([ A(t_1), \psi (t_1 ) ] ) - C' ([ A(t_2), \psi (t_2 ) ] )\\[2mm]
= & \displaystyle{\int_{t_1}^{t_2} \int_Y ( \| \frac{\partial A(t)}{\partial t }
\| ^2 +  \| \frac{\partial \psi (t)}{\partial t } \| ^2 ) d vol_Y}\\[2mm]
=& \displaystyle{ \int_{t_1}^{t_2} \|\nabla C'([A(t), \psi(t)])\|^2_{L^2(Y)} } 
\end{array}
\]
Note that in this case, the finite energy condition for a  connecting orbit
 is the  same as the finite length  condition for the corresponding
 gradient flow. 
\end{R}
 
\subsection{Transversality for ${\cal M}(a, b)$}
\label{transversal}

In this subsection, we will prove that ${\cal M}(a, b)$ is a smooth
manifold after a generic perturbation. As we need $\Bbb R$-translation
action on the moduli space, one may expect  to construct a time-invariant
perturbation of the gradient flow equation to achieve the smoothness for
 ${\cal M}(a, b)$. Unfortunately, this perturbation can't achieve the
transversality for ${\cal M}(a, b)$. In \cite{Fro}, Froyshov construct
an explicit perturbation, here, we show that the 
transversality for ${\cal M}(a, b)$ is actually generic.   

Suppose $a, b$ are irreducible critical points for a fixed perturbation
$\alpha$, choose any smooth representations $(A_1, \psi_1)$, $(A_2, \psi_2)$
in ${\cal A}_{L_1^2}$ for $ a, b$ respectively,
note that ${\cal A}_{L_1^2}$ is the  $L_1^2$-configuration
space, then $(A_i, \psi_i)$ satisfies equation (\ref{SW3'})
\ban
&& *F_A =  q(\psi ) + *d \alpha \nonumber\\[2mm]
&& \dir_A (\psi) =0
\nan

Denote by  ${\cal A}_{L^2_1}(a, b)$ the set 
\[
\left\{
(A, \psi ):\ {\Bbb R} \to {\cal A}_{L^2_1}\left| \begin{array}{l}
(1)\ \lim_{t\to -\infty}(A, \psi )\ \hbox{lies in the gauge orbit of 
             $(A_1, \psi_1)$,}\\[2mm]
(2)\ \lim_{t\to +\infty}(A, \psi )\ \hbox{lies in the gauge orbit of  
             $(A_2, \psi_2)$,}\\[2mm]
(3)\ \displaystyle{
    \int_{ -\infty}^{+\infty}(\|\frac{\partial A(t)}{\partial t}\|^2_{L_1^2}
    +\|\frac{\partial \psi(t))}{\partial t}\|^2_{L_1^2})dt < \infty.}
\end{array} \right. \right\}
\]

It is easy to see that the tangent space of ${\cal A}_{L^2_1}(a, b)$ at
$(A(t), \psi(t))$ is 
\[
T_{(A(t), \psi(t)}({\cal A}_{L^2_1}(a, b))
= L_{1, 0}^2({\Bbb R}, \Omega^1(Y, i{\Bbb R}) \oplus \Gamma(S)) \oplus 
\Omega^0 (Y, i{\Bbb R}) \]

Then ${\cal M}(a, b)$ is the moduli space of the gradient flow
eqation (\ref{SW4t'}) defined on  ${\cal A}_{L^2_1}(a, b)$. 
The perturbation we choose is to add a function of $A(t)$ in 
$L^2_{2, 0}({\Bbb R}, \Omega^1(Y, i{\Bbb R}))$ to the curvature equation in (\ref{SW4t'}).

\ba 
&&\frac{\partial A }{\partial t} = *F_A - q (\psi) -*d\alpha_{A(t)} 
\nonumber\\[1mm]
&&\frac{\partial \psi }{\partial t} =  \dir _A (\psi)
\label{tran}
\na
where $\alpha_{A(t)} \in \alpha + L^2_{2, 0}({\Bbb R}, \Omega^1(Y, i{\Bbb R}))$. 
Note that this perturbed gradient flow equation preserves $\Bbb R$-translation, 
that is, for a fixed $\alpha_t = \alpha_{A(t)}$, if $(A(t), \psi(t))$ is a
solution for (\ref{tran}), then $(A(t+s), \psi(t+s)$ ( for $s\in \Bbb R$)
is also a solution for (\ref{tran}).

Define the parametrized moduli space ${\cal M^P} $ as
\[
{\cal M}^P = \{(A, \psi, \alpha_t ) \ | (A, \psi, \alpha_t )
\  \hbox{is a solution of  (\ref{tran})}\}/{\cal G}
\]

For any solution $[A_0(t), \psi_0(t), \alpha_t ] \in {\cal M}^P$, the 
linearisation of (\ref{tran}) is given by 
\[
 L_{1, 0}^2({\Bbb R}, \Omega^1(Y, i{\Bbb R}) \oplus \Gamma(S)) \oplus 
 L^2_{2, 0}({\Bbb R}, \Omega^1(Y, i{\Bbb R})) \to 
 L_{1, 0}^2({\Bbb R}, \Omega^1(Y, i{\Bbb R}) \oplus \Gamma(S)) 
\]
\[
D : \ \left(
\begin{array}{c}
A \\[2mm] \phi  \\[2mm] \beta_t \end{array}\right)
\mapsto  \frac{\partial}{\partial t}\left(\begin{array}{c}
A \\[2mm] \phi \end{array}\right)
- \left(\begin{array}{cc}
*d & -Dq_{\psi_0}\\[2mm]  .\psi_0 &  \dir _{A_0}
\end{array}\right)
\left( \begin{array}{c}
A\\[2mm] \phi \end{array}\right)
+ \left( \begin{array}{c}
\beta _t\\[2mm] 0 \end{array}\right)
\]
where $\beta_t \in  L^2_{2, 0}({\Bbb R}, \Omega^1(Y, i{\Bbb R}))$, we will choose
$\alpha_t  $ in (\ref{tran}) very small such that  
the index and the spectral flow defined in the next subsection are tha same as
the perturbed operators.

From the surjectivity of $D$, we know  that ${\cal M}^P$ 
is a smooth manifold. Now  apply the Sard-Smale theorem to the projection map
from ${\cal M}^P $ to $\alpha +  L^2_{2, 0}({\Bbb R}, \Omega^1(Y, i{\Bbb R}))$,
for a generic perturbation (consisting an open dense set in 
 $\alpha +  L^2_{2, 0}({\Bbb R}, \Omega^1(Y, i{\Bbb R}))$ \ ), to see that 
 ${\cal M}(a, b) $ is a smooth  manifold whose dimension is given by the 
relative Morse index  defined in the next subsection.

\subsection{Spectral flow and relative Morse index for critical pionts}
 
As in subsection \ref{transversal}, the moduli space of the connecting orbits
$ {\cal M} (a, b)$ is a smooth, compact manifold with the
correct dimension given by the index of the deformation complex. 
In this subsection,
we describe a relative Morse index for the critical points in
 ${\cal M}_{Y, \alpha}$.
In general, as in Floer instanton homology, the Morse index defined
by the dimension of the negative space for the Hessian operator at such
critical points is not well-defined, since the Hessian has infinite dimensional
negative space. It is Floer's  insight to introduce  a relative index to
overcome this diffculty.  The dimension of the moduli space of the connecting
orbits is given by  this relative index.
We adopt the same idea to give the index for our critical points.

In \cite{APS},
Atiyah, Patodi and Singer proved that
\ba
Index (\frac{d}{d t} +  T (t)) = \mbox{ \  spectral flow of
$ T (t)$}
\label{APS}
\na
where $T (t)$ is a path of Fredholm operators with
 invertible limits as $ t \to \pm \infty$.

For any  connecting $ [A(t), \psi(t)] \in {\cal M}  (a, b)$, here
$a, b \in {\cal M}^*_{Y, \alpha}$. The linearization on $\cal A$ is given by
\ba
\frac{\partial }{\partial t} - \left(
\begin{array}{cc}
*d & -Dq_{\psi(t)}\\
.\psi & \dir_{A(t)}
\end{array} \right)
\label{index}
\na
with the limits (as $t\to \pm \infty$):
\ba
\left(
\begin{array}{cc}
*d & -Dq_{\psi _{\pm} }\\
.\psi_{\pm} & \dir_{A_{\pm}}
\end{array} \right)
\label{K(t)}
\na
where $(A_{\pm}, \psi_{\pm})$ are smooth representatives for $a, b$.
These are precisely the Hessian operators at $a, b$ and are  invertible
(note that Hessian at the only reducible critical point is also invertible
see Proposition \ref{generic}). Here  we
think of (\ref{K(t)}) as an operator on the $L^2$ tangent bundle, i.e., fixing
the unique  gauge $G(\tau)$, $\tau \in  \Omega^0(Y, i{\Bbb R})$ given by
(\ref{deform}), such that $K(t) + G(\tau)$ is well defined as a
$L^2$-tangent vector of ${\cal B}^*$. Denote this gauged operator
as $\hat K(t)$. At the critical point, such a $\tau $
is zero. We abuse the notation $K(t)$ and  $\hat K(t)$ when it doesn't cause
any confusion.

By the Atiyah-Patodi-Singer index theorem (\ref{APS}),
the dimension of ${\cal M}  (a, b)$ is determined by
the spectral flow of $\hat K(t)$,  the number of the
eigenvalues crossing  $0$ from negative to positive
 minus the number of the eigenvalues crossing  $0$ from positive to negative.
 In \cite{Wang2}, we use the results of this paper to give a spectral 
flow version of the definition of Casson type invariant and 
 its ${\Bbb Z}_2$-version.
In \cite{Mar2}, M. Marcolli also gave a description of this invariant
for other 3-manifolds.

 To give a $\Bbb Z$-valued index for $a \in {\cal M}_{Y, \alpha}$,
we need to consider a path $\gamma (t)$ from the trivial point $[0, 0]$ to $a$,
then define  the index at $a$ as follows,
\[
\mu (a) =\ \hbox{spectral flow of $\hat K$ along $\gamma (t)$}
\]

There is no ambiguity in the definition of the index $ \mu (a)$, since
for our homology 3-sphere $Y$, $H^1({\cal A}, {\Bbb Z}) =0$. One can see
this by a simple observation dim${\cal M} (a, a)$ is
zero, actually, ${\cal M}  (a, a)$ is just one point,
the constant solution, by the identity
\ban 
&&\quad\int^{+\infty}_{-\infty} \|\nabla C'([A(t), \psi(t)])\|^2_{L^2(Y)}\\[2mm]
&&= \|\displaystyle{\frac{\partial A(t)}{\partial t}}\|^2_{Y\times {\Bbb R}} 
+ \|\displaystyle{\frac{\partial \psi(t))}{\partial t}}\|^2_{Y\times {\Bbb R}} \\[2mm]
&&= \lim_{t \to \infty}C'(A(t), \psi(t)) - \lim_{t \to -\infty}C'(A(t),\psi(t))
\nan
which is zero for the solution $ [A(t), \psi(t)]$ 
in ${\cal M} (a, a)$.

This will make the Seiberg-Witten-Floer homology graded by 
  $\Bbb Z$. From the Atiyah-Patodi-Singer index theorem, we have the follow lemma.

\begin{Le}
For $a, b \in {\cal M}^*_{Y\ \alpha}$ with $\mu(a) > \mu(b)$, the moduli
space of the ``connecting orbits" ${\cal M} (a, b)$
 is a $\mu(a)-\mu(b)$ dimensional,
 smooth, compact  manifold. For $a \ne  b$ with  $\mu(a) \le  \mu(b)$,
 ${\cal M}  (a, b)$ is empty. 
\end{Le}
\begin{proof} It is sufficient to prove
 that for $a \ne b$ with  $\mu(a) =\mu(b)$
 then ${\cal M} (a, b)$ is empty, this follows
 from  the equations 
(\ref{tran}) being invariant under the action of the $\Bbb R$-translation, if
 ${\cal M}  (a, b)$ is non-empty, its dimension has to
be at least 1, this is impossible from the index calculation.
\end{proof}

\begin{R}
We didn't discuss the dimension of ${\cal M} (a, b)$
where either $a$ or $b$ is the only reducible critical point (which is 
also non-degenerate according to Lemma \ref{iso}). By the excision 
principle, 
\[
dim {\cal M}  (a, b) = \mu(a) -\mu(b) -1
\]
\[
dim {\cal M} (c, a) =  \mu(c) -\mu(a) 
\] 
where $a$ is the reducible critical point. 
We delay this proof until we introduce our  gluing theorem. 
\label{reddim}
\end{R}

\subsection{Decay estimates for the gradient flow near the critical points}

In this subsection, we will study the gradient flow  near the critical 
point $[A_0, \psi_0] $ of $C'$, $\psi_0 \ne 0$,
 where $[A_0, \psi_0] $ obeys
\ba
&& *F_A =  q(\psi ) + *d \alpha \nonumber\\[2mm]
&& \dir_A (\psi) =0
\label{A}
\na
By Lemma \ref{TB} (2), we know that $U([A_0, \psi_0], \epsilon)$
is diffeomorphic to a neighbourhood $U(0, \epsilon)$
 of $T_{[A_0, \psi_0]}({\cal B}^*)$. Using this
relation, we can rewrite the Seiberg-Witten equations (\ref{SW4t'})    
on $T_{[A_0, \psi_0]}({\cal B}^*)$ as follows,
\[ 
\frac{\partial}{\partial t}\left(\begin{array}{c}
A_0 + A \\[2mm] \psi_0 + \phi \end{array}\right)
=\left(\begin{array}{c}*F_{(A_0 +A) } - q(\psi_0 + \phi )- *d\alpha \\[2mm]
\dir_{(A_0 + A)}(\psi_0 + \phi ) \end{array}\right)
\]
Applying (\ref{A}), we simplify the gradient flow equation on
 $T_{[A_0, \psi_0]}({\cal B}^*)$ as follows:
\ba
\frac{\partial}{\partial t}\left(\begin{array}{c}
A \\[2mm] \phi \end{array}\right)
=\left(\begin{array}{cc} 
*d & -Dq_{\psi_0}\\[2mm]  .\psi_0 &  \dir _{A_0}
\end{array}\right)
\left( \begin{array}{c}
A\\[2mm] \phi \end{array}\right) 
+ \left( \begin{array}{c}
-q(\phi) \\[2mm] A.\phi  \end{array}\right)
\label{ODE}
\na

Write
the quadratic term in (\ref{ODE}) as $Q(A, \phi)$. 
Denote 
\[
\left(\begin{array}{cc} 
*d & -Dq_{\psi_0} \\[2mm]  .\psi_0 &  \dir _{A_0}
\end{array}\right)
\] 
as $K_{[A_0, \psi_0]}$, it is the Hessian operator of $C'$ at $[A_0, \psi_0]$
the linearization of the perturbed Seiberg-Witten equations on $Y$. 
$K_{[A_0, \psi_0]}$ is a closed, self-adjoint, Fredholm operator from the
$L^2_1$-completion of $ T_{[A, M]}({\cal B})$ to
 the $L^2$-completion of $ T_{[A, M]}({\cal B})$, and  has only discrete spectrum
without accumulation points.

From the preceding sections, we know that $K_{[A_0, \psi_0]}$ is 
a invertible operator acting on $T_{[A_0, \psi_0]}({\cal B}^*)$, we can decompose
$T_{[A_0, \psi_0]}({\cal B}^*)$ as the direct sum of the positive
eigenvalue spaces and the negetive eigenvalue space. 

\ba
T_{[A_0, \psi_0]}({\cal B}^*) =( \oplus _{\lambda >0} {\cal H} _{\lambda })
\oplus (\oplus _{\lambda < 0} {\cal H} _{\lambda }) 
\label{dec}
\na
where ${\cal H} _{\lambda }$ is the eigenspace with eigenvalue $\lambda $.

Fix $\delta  < \min \{ |\lambda |\}$, let the semi-group generated
by $K_{[A_0, \psi_0]}$ be
$\Phi (t_0,  t)$, then $\Phi (t_0, t)$ preserves the decomposition (\ref{dec})
We  write $ x \in T_{[A_0, \psi_0]}({\cal B}^*)$
as 
\[ x = x^+ + x^- \]
where $x^+ \in  {\cal H}^+=\oplus _{\lambda >0} {\cal H} _{\lambda }$ is called
the stable part, $  x^- \in {\cal H}^-= \oplus _{\lambda< 0} {\cal H} _{\lambda }$ is
called the unstable part. Denote ${\cal H}  = {\cal H}^+ \oplus {\cal H}^-$.

We have the following hyperbolicity formula for $\Phi (t_0,  t)$:

\begin{itemize}
\item  For  $x^+ \in  {\cal H}^+, t \ge t_0$, we have
\[
\| \Phi (t, t_0) x^+ \| \le  \exp^{-\delta (t-t_0)} \|x^+\|.
\]
\item  For  $  x^- \in {\cal H}^-, t \le t_0$, we have \[
\| \Phi (t, t_0) x^- \| \le  \exp^{-\delta (t_0 -t)} \|x^- \|.
\]
\end{itemize}

We can solve the equation (\ref{ODE}) near  $0$ in ${\cal H}^+ \oplus {\cal H}^-$,
let 
\[ 
x = \left(\begin{array}{c}
  A \\ \phi \end{array}\right) \]
satisfy the following equation with some boundary data:
\ba
\left\{\begin{array}{l}
       \displaystyle{\frac{\partial x}{\partial t}} 
                = K_{[A_0, \psi_0]} x + Q(x)\\[5mm]
       \| Q(x) \|_{L^2} \le C_1\|x\|^2_{L^2_1} \\[5mm]  
       \|D_x Q(x)\|_{L^2} \le \epsilon \\[5mm]
       Q(0) =0.
       \end{array}
\right. 
\label{hyper}
\na

Generally, the solution is described by the following integral equation,
\[
x(t) = \Phi (t_0,  t)x(t_0) + \int^t_{t_0}\Phi (s,t)Q(x(s))ds.
 \]

The following two lemmas can be obtained by  applying the theory developed in
\cite{AB}, which is a nice model for decay estimates.

\begin{Le}
For sufficiently small $p, q$ in ${\cal H}^+, {\cal H}^-$ respectively, 
$T_1 < T_2 $, there is a unique solution $x(t): \ [ T_1, T_2]
\to {\cal H}$ satisfying
(\ref{ODE}) such that 
\[\left\{\begin{array}{l}
x^+ (T_1) =p\\[3mm]
x^-(T_2) = q \\[3mm]
\| x(t) \|_{L_1^2} < \epsilon.
 \end{array}
\right. \]
This solution depends smoothly on $p, q, T_1,  T_2$. Denote
$x(t) = x(t, p, q, T_1, T_2)$. 
\label{unique}
\end{Le}
\begin{proof}
The solution must have the following form
\ba
x(t) &=&  \Phi (T_1, t)p + \int^t_{T_1}\Phi (s,t)Q^+(x(s))ds\nonumber \\[2mm]
     && + \Phi (T_2, t) q - \int^{T_2}_t \Phi (s,t)Q^-(x(s))ds.
 \label{sol}
\na
Write (\ref{sol}) as a fixed point equation:
\[
x = F(x)\]
where $x$ lies in the Banach space $C^0([T_1, T_2],  {\cal H})$ 
with the supremum norm. 
The existence and uniqueness can be obtained by showing that the above map
is a strong contraction on a sufficiently small ball $\|x\| < \eta$. This can be verified by  (\ref{hyper}) as follows,
\[\begin{array}{rcl}
\| F(x)\|& <& C(\displaystyle{ \epsilon + |\int_{T_1}^{T_2}\exp^{\delta s}ds|\sup_{\|x\| < \eta}
\|Q(x)\|}) \\[2mm]
&<&C(\epsilon +\displaystyle{\frac{C_1}{\delta}}\eta ^2) < \frac{1}{2}\eta
\end{array}
\]
provided that $\epsilon$ is sufficiently small. 
\end{proof}

\begin{Le}
Suppose $x(t)$ is the solution in Lemma \ref{unique}, then there are 
a priori estimates for $x(t)$ as follows where $C_3$ is a positive constant,
\[   \begin{array}{l}
\| x^+(t) \|_{L_1^2}  \le  C_3  \exp^{-\delta (t-T_1)}\\[3mm]
\| x^-(t) \|_{L_1^2}  \le C_3  \exp^{-\delta (T_2 -t )}\\[3mm]
\| \displaystyle{\frac{\partial x^+(t)}{\partial p}} \|_{L_1^2} 
 \le  C_3\exp^{-\delta (t-T_1)}\\[3mm]
\| \displaystyle{\frac{\partial x^+(t)}{\partial q} } \|_{L_1^2}  \le
 C_3 \exp^{-\delta (t-T_1)}\\[3mm]
\| \displaystyle{\frac{\partial x^-(t)}{\partial p}} \|_{L_1^2} 
 \le C_3  \exp^{-\delta (T_2-t )}\\[3mm]
\| \displaystyle{\frac{\partial x^-(t)}{\partial q}	} \|_{L_1^2} 
 \le  C_3\exp^{-\delta (T_2-t )}
  \end{array}
\]
\label{est}
\end{Le}

Note that from the first two estimates, it is easy to see that 
the solution of (\ref{ODE}) satisfies the following property,
\[
\| x(t)\|_{L_1^2}\le C_3 \exp^{-\delta d(t)}
\]
where $d(t) = \min \{t-T_1, T_2 -t\}$.

\begin{proof}
These estimates can be obtained by the standard open and closed
argument and continuity. We only prove the first estimate.
 Define the non-empty open set 
\[
S=\{t \in [T_1, T_2]| \| x^+(t) \|_{L_1^2}  \le  C_3  \exp^{-\delta (t-T_1)}
\]
The aim is to prove $S$ is also closed.
Since $T_1 \in S$, suppose $ [T_1, t ) \subset S$.
 Note that $ x^+(t)$ can be write
as 
\[
x^+(t)= \Phi (T_1, t)p +\int_{T_1}^t\Phi (s, t)Q^+(x)ds
\]
Estimating the norm  $x^+(t)$ 
shows that $t\in S$, hence $S =[T_1, T_2]$.
\end{proof}

\begin{Pro}
There are positive constants $\epsilon, \delta, C_4$, such that for any $
T_2 > > T_1 $, if $[A(t), \psi (t)] $ is the solution for the perturbed
Seiberg-Witten equations in a temporal gauge on $[T_1, T_2] \times Y$ near
the critical point $[A_0, \psi_0]$, that is, if for each $t \in [T_1, T_2]$
\[
dist_{L^2_1}([A(t), \psi (t)], [A_0, \psi_0]) < \epsilon
\]
then there is a exponential decay estimate as follows,
\[
dist_{L^2_1}([A(t), \psi (t)], [A_0, \psi_0]) \le C_4  \exp^{-\delta d(t)}
 \]
where  $d(t) = \min \{t-T_1, T_2 -t\}$.
Therefore, if the gradient flow $[A(t), \psi (t)]|_{[T, \infty )}$ is 
sufficiently close to the 
critical point $[A_0, \psi_0]$, then  as $ t \to \infty $, $[A(t), \psi (t)]$
decays to  $[A_0, \psi_0]$ exponentially. There is a similar exponential
decay for the finite energy solution $[A(t), \psi (t)]$ as $t \to -\infty $. 
\label{decay}
\end{Pro}

\begin{proof} 
We have  identifed 
 the perturbed Seiberg-Witten equations in temporal gauge and
(\ref{ODE}) on $\cal H$,
then this proposition is a direct consequence of this identification and
the above lemmas for the small solutions.
\end{proof}

The following lemma gives a nice picture for the limits of the solution near the
critical  point, the broken trajectories appear in the limit. This will play
an important role in the definitions of the boundary operators for
the Floer
complex and various chain-maps for metric-independence.

\begin{Le}
Let $t\mapsto x(t, p, q, -T, T)$ be the solution as in Lemma \ref{unique},
As $T \to \infty $,
\\
(1) The trajectories: $\gamma_1 (t): \ [0, T] \to {\cal H}$
\[
\gamma_1(t) = x (t-T,  p, q, -T, T)\]
 approach the limit which lies in 
${\cal H}^+$ (locally stable manifold).
\\
(2) The  trajectories: $\gamma_2 (t): \ [0, T] \to {\cal H}$
\[
\gamma_1(t) = x (T-t,  p, q, -T, T)\]
 approach the limit which lies in
${\cal H}^-$ (locally unstable manifold).
\label{broken}
\end{Le}

\begin{proof}
From Lemma \ref{est}, we know that 
\[
\| \gamma_1(t)^-\|_{L^2_1} = \|x^- (t-T,  p, q, -T, T)\|_{L^2_1} 
\le  C_3  \exp^{-\delta (T-(t-T))}
\]
\[
\| \gamma_2(t)^+\|_{L^2_1} = \|x^+ (T-t,  p, q, -T, T)\|_{L^2_1}
\le  C_3  \exp^{-\delta (T-t-(-t))}
\]
then this lemma follows.
\end{proof}

\subsection{Gluing arguments}

In this subsection, we construct the gluing map which will enable us to 
build  a Floer complex and various chain homomorphisms. For simplicity,
we denote 
${\cal M}_{Y, \alpha}$ as ${\cal R}$. 

Since there is a natural $\Bbb R$-action on  ${\cal M} (a, b)$ (time translation), 
we define the $\Bbb R$-quotient space by,
 \[ \hat {\cal M} (a, b) =  {\cal M} (a, b)/ {\Bbb R} .\] 

\begin{Pro}
For  $T$ sufficiently large, 
let $a, b, c \in {\cal R}$ with $\mu(a) > \mu (b) > \mu(c)$. Suppose
$b$ is irreducible, then there is an embedding map, for $T$ sufficiently large,
\[
g : \qquad \hat {\cal M} (a, b) \times \hat { \cal M} (b, c) \times
  [ T, \infty )\to  \hat {\cal M} (a, c). \]
 If $b$ is the reducible critical point, then there is a gluing group $U(1)$
arising from the stabiliser of $b$ in the gluing map; this means
that there is a local
diffeomorphism,
\[
\hat g: \qquad \hat {\cal M} (a, b) \times  \hat {\cal M} (b, c) \times 
U(1)\times [ T, \infty )\to  \hat {\cal M} (a, c).
\]
\label{glue}
\end{Pro}

\begin{proof}
Choose $x(t) \in \hat{\cal M} (a, b)$ and $ y(t) \in \hat{\cal M} (b, c)$, 
by Proposition \ref{decay}, assume $b=[A_0, \psi_0]$,
 we can write (for T sufficiently large) 
\[
x(t)= b + \left(\begin{array}{c}
          A_1 \\ \phi_1 \end{array}\right) \]
\[
y(-t) = b + \left(\begin{array}{c}
          A_2 \\ \phi_2 \end{array}\right) \]
for $t > 2T$, then $ x_i(t)= (A_i, \phi_i ) $ satisfies (\ref{ODE}), that is,
\[
\frac{\partial}{\partial t}\left(\begin{array}{c}
A_i \\ \phi_i \end{array}\right)
=\left(\begin{array}{cc}
*d & -Dq_{\psi_0}\\  .\psi_0 &  \dir _{A_0}
\end{array}\right)
\left( \begin{array}{c}
A_i\\ \phi_i \end{array}\right)
+ \left( \begin{array}{c}
-q(\phi_i) \\ A_i.\phi_i  \end{array}\right).
\]
We simplify this equation as 
\ba
 \frac{\partial x_i}{\partial t} = Kx_i + Q(x_i) 
\label{ODE1}
\na
with the exponential decay estimates,
\[
\|( A_i, \phi_i )\|_{L^2_1} < C_5 \exp^{-\delta T}. \]
We construct the following gluing path on $\cal B$,
\[
x\sharp_Ty(t)=
\left\{\begin{array}{l@{\qquad \hbox{for}\qquad }c}
x(t+2T) & t\le -1 \\[2mm]
b+\rho(t) x_1(t+2T) + (1-\rho(t))x_2(t-2T) & -1\le t \le 1 \\[2mm]
y(t-2T) & t\ge 1
\end{array}
\right. \]
where the calculation is done near $b$ under the identification 
$U(b, \epsilon)$ with a domain in $\cal H$, $\rho(t)$ is the cut-off
function on $[-1, 1]$,
whose derivative has support in $[-\frac 12, \frac 12]$, is $1$ 
near $-1$ and 0 near 1. Assume $|\partial \rho| < C_6$.
We will find a solution for (\ref{SW4t'}) uniquely determined by the 
gluing data $x(t), y(t), T$.  

For $t\in [-1, 1]$, write  
\[
x'(t) = x\sharp_Ty(t) -b = \left( \begin{array}{c}
        A' \\ \psi'
\end{array}\right)
\]
then for $T$  sufficiently large, $  x\sharp_Ty(t) -b $ is an approximate
solution for (\ref{ODE1}), in the sense that,
$$ \eta = \partial _t x' -Kx' -Q(x') $$ 
is sufficiently small as long as $T$ is sufficiently large.

We aim to find a $\xi'(t) \in L^2_{1, 0}( [-1, 1],  {\cal H}) $ (very small) such
that  $  x' + \xi'$ solves the flow equation. This means that $ \xi'$
obeys the following equation,
\[
-\frac{\partial \xi'}{\partial t} + K\xi' + Q(\xi') 
+ \left( \begin{array}{cc}
0 &-Dq_{x'}   \\[1mm]
. \psi' & A'
      \end{array} \right) \xi' 
= \eta .
\]

Note that, for T sufficiently large: 
\begin{itemize}
\item $ S =\left( \begin{array}{cc}
0 &-Dq_{x'}  \\[1mm]
. \psi' & A'
      \end{array} \right) $ has arbitrarily  small operator
norm on $ L^2_{1, 0}( [-1, 1],  {\cal H}) $.  

\item $-\displaystyle{\frac{\partial }{\partial t}} + K $
 is a bounded, invertible 
operator on the following space,
\[
 L^2_{1, 0}( [-1, 1],  {\cal H}) \longrightarrow 
 L^2_0( [-1, 1],  {\cal H}) .
\]
\end{itemize}

Only the second assertion needs some explanation. Since $K$ is invertible,
if  $$z(t) \in L^2_{1, 0}( [-1, 1],  {\cal H})$$ in 
$ker (-\displaystyle{\frac{\partial }
{\partial t}} + K) $, then $z(t)$ satisfies
\[
\left\{\begin{array}{l}
\displaystyle{\frac{\partial z(t) }{\partial t}} = K z(t)\\[2mm]
z(0)=0\\[2mm]
z(1) =0 \end{array}
\right. \]

From the flow $\Phi (0, t)$ generated by $K$, we see that $z(t)$ must be zero.
Then $ -\displaystyle{\frac{\partial }{\partial t}} + K $ is an  invertible  
operator. The boundness is obvious. So is 
 $-\displaystyle{\frac{\partial }{\partial t}} + K + S$. 
 Therefore, there is a right inverse 
$P$, $P$ is a bounded operator acting on $ L^2_0( [-1, 1],  {\cal H})$ 
with range $ L^2_{1, 0}( [-1, 1],  {\cal H}) $.
  
We write $ \xi' = P(\xi )$, then $\xi$ satisfies 
   \[ \xi + Q(P(\xi)) = \eta \]
where $\eta $ can be made arbitrarily  small by making $T$ larger.

Since $Q$ is a quadratic function, then we have 
\[
\| Q(P(\xi_1))-Q(P(\xi_2))\| \le C_7 \|\xi_1 -\xi_2\| (\|\xi_1\| 
+\|\xi_2\|)\] 
for some constant $C_7$. Then the existence and uniqueness of the solution
determined by $x(t), y(t), T$ follows from the map
\[ \xi \mapsto \eta - P(Q(\xi ))\]
is a strong contraction. The fixed point of this map,  $\xi$,
is our $g (x(t), y(t), T)$. $g$ is an injection by our construction. 
It is easy to see that the image of $g$ is one end of $ \hat {\cal M}$.

When gluing at the reducible solution, basically the same procedure applies.
We need to work with the framed  configuration space near the reducible
critical point, 
the stabiliser of the reducible  solution then enters the gluing map by a
free action. We only construct the $\sharp$-map in this case, gluing 
two gradient flows $x(t), y(t)$ along a reducible point $b$,
 suppose $ T$ is sufficeintly large, $\exp (i\theta )\in U(1)$, then as before
\[
x\sharp_T^\theta y(t)=
\left\{\begin{array}{l@{\quad \hbox{for}\quad }c}
x(t+2T) & t\le -2 \\[2mm]
b+ \exp^{i\theta(t+2)}.x_1(t+2T) & -2 \le t \le -1\\[2mm]
b+\exp^{i\theta}. (\rho(t) x_1(t+2T) +
 (1-\rho(t))x_2(t-2T)) & -1\le t \le 1 \\[2mm]
b+exp^{i\theta (t-2)}.x_2(t-2T) & 1\le t\le 2\\[2mm]
y(t-2T) & t\ge 2
\end{array}
\right. \]
It is easy to see that $x\sharp_T^\theta y(t)$ is an approximate solution, to 
find the unique correct solution defined by $x, y, T, \theta$, one follows
 the steps
we did for the irreducible gluing map. The details we leave to the reader.
\end{proof}

This proposition gives an explicit picture for the boundary of   various
moduli spaces, we give here just one example of these arguments, we will meet 
several other interesting moduli spaces when we discuss  
  relative Seiberg-Witten invariants. Proposition \ref{glue}
also tells us how a  trajectory connecting two critical points $a, b$ 
($\mu (a) > \mu (b)$) breaks into two pieces, 
which become  two
 trajectories satisfying (\ref{SW4t'}) and breaking at another critical point
$c$ with $\mu(a) < \mu (c) < \mu (b)$.

 From this Proposition, we can prove 
\[ dim {\cal M}  (a, b) = \mu(a) -\mu(b) -1
 \]
and \[  dim {\cal M}  (c, a) = \mu(c) -\mu(a) \]
where $a$ reducible.
Suppose ${\cal M}  (c, a), {\cal M} (a, b) $ are  nonempty moduli spaces, then
the gluing map $\hat g$ in Proposition \ref{glue} tells us the identity,
\[
dim {\cal M}  (c, a) -1 + dim {\cal M}  (a, b)-1 +2 = dim {\cal M}  (c, b)-1 \]
since $c, b$ must be irreducible, we have $ dim {\cal M}  (c, b)=\mu(c) -\mu(b)$, the assertion
follows. In particular, when $\mu(a) -\mu(b) =1$ and $a$ is reducible,
$ {\cal M}  (a, b)$  is empty.    

\begin{Cor}
Suppose $ a, c$ are two irreducible critical points with the  index given by 
$\mu(a) = \mu(c)+2$,
the boundary of $ \hat {\cal M}(a, c)$ consists of the union 
\[
\bigcup_{\{b| \mu(b)= \mu(a) -1\}}\hat{\cal M}(a, b) \times \hat{\cal M}(b, c)
\]
where $b$ runs over  the set of the irreducible critical point.
\label{boundary}
\end{Cor} 

\begin{proof}
From the gluing map (see Proposition \ref{glue})and the above arguments,
we know that  the limit of a sequence of trajectories
doesn't break through the reducible critical point, and only breaks at the
irreducible critical point $b$ with $\mu(b)= \mu(a) -1$.
\end{proof}

\section{Floer Homology}

We use the simplified notations as before, 
such as ${\cal R}, {\cal R}^*, {\cal M}(a, b), \hat{\cal M}(a, b)$, where
$ {\cal R}^*$ is a finite set of irreducible, non-degenerate critical
points, indexed by the spectral flow of a path starting from the trivial point.
For $a, b \in  {\cal R}^*$ with $\mu(a) -\mu(b) =1$, then $\hat{\cal M}(a, b)$
is also a finite set of points, specifically,
we recall $\hat{\cal M}(a, b)$ defined by 
 \ba
 \hat{\cal M}(a, b)  = \left\{
x(t)\left| \begin{array}{l}
(1) \hbox{$x(t)$ is the finite energy solution of (\ref{SW3'})} \\[2mm]
    \quad \hbox{ modulo the gauge group and $\Bbb R$-translation,}\\[2mm]
(2) \lim_{t \to +\infty}x(t)= a \in {\cal R}^*,\\[2mm]
(3) \lim_{t \to -\infty}x(t)= b  \in {\cal R}^*.
\end{array} \right. \right\}
\label{hat}
\na

\begin{The}{\bf Seiberg-Witten-Floer homology:} 
Let $C_k(Y) $ be the free abelian group over $\Bbb Z$ with generators consisting
of points in $ {\cal R}^*$ whose index is $k$.  
Let $n_{a\ b}$ be the sum of $\pm 1$ over $\hat{\cal M}(a, b)$ whenever
$\mu(a)-\mu(b) =1$, where the sign is determined by comparing the orientation
on $ {\cal M}(a, b)$ with the $\Bbb R$-translation on each isolated 1-dimensional
component, the orientation on $ {\cal M}(a, b)$ is determined by the 
orientation of $Y$. 
Then the boundary operator 
\[
\partial : \qquad C_k(Y) \longrightarrow C_{k-1}(Y) \]
\[ \partial (a) = \sum_{b\in {\cal R}^*\ \mu(b) = k-1} n_{a\ b}b \]
satisfies $\partial\partial =0$. The homology group 
\[
HF^{SW}_k(Y,g) = ker \partial _k/Im \partial_{k+1} \]
is a topological invariant (up to a canonical isomorphism). Moreover, 
$HF^{SW}_*$ is a functor on the category of homology 3-spheres and their 
cobordisms.
\label{HF}
\end{The}

\begin{R}
The Casson-type invariant for a homology 3-sphere  defined in 
\cite{Wang2} is the Euler characteristic of the Seiberg-Witten-Floer
homology $HF^{SW}_*$:
\ba
\begin{array}{rcl}
\chi (HF^{SW}_*(Y, g) ) &=& \sum _k (-1)^k \dim HF^{SW}_k(Y, g)\\[2mm]
&=& \sum_k  (-1)^k  \dim C_k (Y)\\[2mm]
&=& \sum_{a\in {\cal R}^*} (-1)^{\mu(a)}\\[2mm]
&=& \lambda (Y, g)
\end{array}
\label{casson}
\na
A similar formula was found by M. Marcolli \cite{Mar2} for 
a 3-manifold  with $Spin^c$ structure whose determinant line bundle has
non-zero first Chern class.
\end{R}

\begin{R}
Obviously, for a diffeomorphism $f: Y \to Y$, there is an isomorphism
\[
HF^{SW}_*(Y, g) \cong HF^{SW}_*(Y, f^*g)\]
\end{R}
 
\begin{proof}
We only prove $\partial\partial =0$, the remaining assertions will be proved in 
the following subsection.
By definition, 
\[ \begin{array}{rcl}
\partial^2 (a)&=& \displaystyle{\sum_{b\in {\cal R}^*\ \mu(b)=\mu(a) -1}}
 n_{a \ b} \partial(b)\\[2mm]
&=& \displaystyle{\sum_{b\in {\cal R}^*\ \mu(b)=\mu(a) -1}\quad
 \sum_{c\in {\cal R}^*
 \ \mu(c)= \mu(b)-1}}n_{a \ b}n_{b \ c}c
\end{array}
\]
 To see $$ \displaystyle{ \sum_{b\in {\cal R}^*\ \mu(b)=\mu(a) -1}
 n_{a \ b}n_{b \ c}=0}$$ for any
$a, c \in {\cal R}^*$ whenever $\mu(a)= \mu(c)+2$. 
We know that  the number $$\displaystyle{\sum_{b\in {\cal R}^*\ \mu(b)=\mu(a) -1}
n_{a \ b}n_{b \ c}}$$ 
is the number of oriented boundary points of $ \hat {\cal M} (a, c)$
(Proposition \ref{boundary} ), hence is zero.
Now the Seiberg-Witten-Floer homology group is well-defined.
\end{proof}
  
To understand the metric (perturbation) dependence of $HF^{SW}_k(Y, g)$,
we should consider the equivariant Seiberg-Witten-Floer homology
$HF^{SW}_{k, U(1)} (Y, g) $ (see \cite{MW} for the definition), which
is independent of metric and perturbation up to index-shifting. Here we 
apply the results in \cite{MW} to illustrate how the irreducible
critical points interact with the reducible
one (denoted by  $\theta$). 

\begin{The} (Proposition 7.8, Proposition 8.1 in \cite{MW}) \\
For $k< 0 $,
\[ HF^{SW}_k (Y, g_0)  \cong HF^{SW}_{k, U(1)}(Y, g_0). \]
For $k\ge 0$, we have the following exact sequences
\[
0\to HF^{SW}_{2k+1, U(1)}(Y, g_0) \stackrel{i_{2k+1}}{\to}
 HF^{SW}_{2k+1}(Y, g_0) 
\stackrel{\bigtriangleup_{k}}{\to}  {\Bbb R }\Omega^k
\to \]\[
\to HF^{SW}_{2k, U(1)}(Y, g_0) \to HF^{SW}_{2k}(Y, g_0)\to 0.
\]  
where $\bigtriangleup_{k}$ is given by 
 \[ \bigtriangleup_{k}(\sum_a x_a a) =x_a m_{ac} m_{ce}\cdots
m_{\alpha' \alpha} n_{\alpha\theta} \Omega^k
\]
Here  a sum over repeated indices is  over all critical
points with indices $\mu(a)=2k+1$, $\mu(c)=2k-1$,
$\mu(\alpha')=3$, $\mu(\alpha)=1$ and $m_{\alpha \beta}$ for
critical points $\alpha, \beta$ with relative index $2$ is the Seiberg-Witten
invariant on $Y \times {\Bbb R}$ with boundary conditions $\alpha, \beta$
respectively. 
\end{The}

Let $\alpha$ be a critical point with index $2k+1$, donote 
\[
 m(\alpha, \theta ) = 
m_{\alpha a} m_{ac} m_{ce}\cdots m_{\beta' \beta}n_{\beta, \theta }
\]
the number presented in $\bigtriangleup_{k}$, sum over $\mu(a)=2k-1$,
 $\mu(c)=2k-3$, $\cdots$, $\mu(\beta') =3$, $\mu(\beta) = 1$.
Then we see that 
\[
\begin{array}{lll}
Im (i_{2k+1}) & = & Ker (\bigtriangleup_{k})\\[2mm]
&=& \left\{ \sum_{\alpha } x_\alpha  \alpha  \left| \begin{array}{ll}
(1) & \ \sum_{\alpha } x_\alpha n_{\alpha \ \beta } =0, \hbox{ for $\beta,
\mu(\beta) = 2k-1$} ,  \\[2mm]
(2) & \ \sum _{\alpha } x_\alpha m(\alpha, \theta ) =0 \end{array} 
\right. \right\}
\end{array} \]

Therefore, $Ker  (\bigtriangleup_{k})$ measure the interaction of
$HF^{SW}_*(Y, g)$ with the reducible critical points.  
There is a similar analogue for Seiberg-Witten-Floer cohomology.
  
\section{Relative Seiberg-Witten invariants}

In this section, we define the relative Seiberg-Witten invariant
for a 4-manifold with a cylindrical end (a homology 3-sphere), this invariant 
takes its values in the Seiberg-Witten-Floer homology group we defined in 
Theorem \ref{HF}. In instanton theory, this is true for a general
manifold see \cite{D1} \cite{BD}. From the result of Marcolli,
there is also a well-defined primary  relative Seiberg-Witten invariant valued
in the  Seiberg-Witten-Floer homology group for 3-manifold 
with $Spin^c$ structure
whose determinant line bundle has
non-zero first Chern class.

 First we give the definition of  
the finite energy solution on a 4-manifold with a cylindrical end
which is isometric to $[0, \infty ) \times Y$ for a homology 3-sphere $Y$.
Where the boundary is $S^1 \times \Sigma_g \ (g>1)$ for a Riemnanian 
surface $ \Sigma_g$, this definition appears in \cite{MST}.

\begin{D}
Let $X$ be a 4-manifold with end isometric to $[0, \infty )
\times Y$ for a homology 3-sphere $Y$. 
Fix a $Spin^c$ structure on $X$ whose restriction to the cylindrical end
$[0, \infty ) \times Y$ is the pull-back  of the $Spin^c$ structure on $Y$ with
the line bundle we introduced earlier. For any 
solution $(A, \psi )$ to the Seiberg-Witten equations  on $X$
 with respect to this
$Spin^c$ structure, in a temporal gauge on the cylindrical
end, there is  a unique 
gradient flow $x(t): [1, \infty ) \to {\cal B}^*$, determined by 
the solution $(A, \psi )$. A
  finite energy solution to the Seiberg-Witten equations on $X$ is a 
solution for which the associated gradient flow line 
 $x(t): [0, \infty ) \to {\cal B}^*$ 
satisfies the condition that 
\[
\lim_{t\to \infty} (C'(x(t)) - C'(x(1)))
\]
is finite.
\label{cylin}
\end{D}

In the above definition, the (perturbed) Seiberg-Witten equations on $X$ are 
\ba 
&&F^+_A = \frac{1}{4}<e_ie_j\psi, \psi>e^i\wedge e^j + \rho (d\alpha
+ dt\wedge *d\alpha )\nonumber\\[2mm]
 &&\dirac_A(\psi ) =0
\label{SWc}\na
where  $*$ is the complex Hodge star operator on $Y$,  the cut-off 
function $\rho $ has support in $(0, \infty ) $, and equals   $1$ over 
$[1, \infty )$. 

Define ${\cal M} (X, a)$ as the moduli space with asymptotic limit $a\in {\cal R}$
\ba
{\cal M} (X, a)= \left\{
[A, \psi ]\left| \begin{array}{l}
(1)\ [A, \psi] \ \hbox{ denotes the gauge orbit of }\\[2mm]
\qquad \hbox{ the solution $(A, \psi)$
 of (\ref{SWc})}
\\[2mm]
(2)\ \hbox{$x(t)$ is the finite energy solution on $Y \times [1, \infty )$}
\\[2mm]
\qquad  \hbox{ associated with $(A, \psi)$ as in Definition \ref{cylin}}
\\[2mm]
(3)\ \lim_{t \to\infty}x(t)= a 
\end{array} \right. \right\}
\label{Xa}
\na

Similar to the exponential decay in Proposition \ref{decay}, any finite
energy solution to the Seiberg-Witten equations (\ref{SWc}) 
decays exponentially to a critical point in $\cal R$, from this fact we can 
get a similar gluing theorem from which we draw the following conclusion.

\begin{Pro}
After a small compact support perturbation of (\ref{SWc}), 
${\cal M}(X, a)$ is a smooth manifold, whose
boundary consists of 
\[
{\cal M}(X, b) \times \hat{\cal M}(b, a)
\] 
where $b$ runs over $\cal R$ with $\mu(b) \ge \mu(a) +1$ and
$\hat{\cal M}(b, a)$ is defined by (\ref{hat}).
\label{boundary4}
\end{Pro}

For any $a \in {\cal R}^*$, we can attach an invariant for the 0-dimensional
component of ${\cal M}(X, a)$, called $n_{X,\ a}$.
Now we define the relative Seiberg-Witten invariant 
with values in the Floer complex,  by
\ba
SW_{X, L}  = \sum_{a\in {\cal R}^*} n_{X, \ a} a
\label{relative}
\na
then from Proposition \ref{boundary4} and Remark \ref{reddim},
 $ SW_{X, L}$ is closed under
the Floer boundary operator, that means
\[
SW_{X, L}  \in HF^{SW}_*(Y).
\]

For a closed, oriented four manifold $X = X_1 \cup_YX_2$, 
which splits along a homology 
3-sphere $Y$, we put a $Spin^c$ structure on $X$ whose restriction to $Y$
is $L$ and then define the Seiberg-Witten-Floer homology $HF^{SW}_*(Y)$.
Then the relevant invariant of 
$X$ is given by the following gluing formula
\[
SW_X = \sum_{a \in {\cal R}^*} n_{X, \ a}n_{X_2, \ a}.
\]

\begin{Ex} As an example, we apply these relative Seiberg-Witten invariants 
to study certain homology spheres which bound  Stein surfaces.
For the background of contact structures and Stein surfaces, see 
\cite{E1} \cite{E2} \cite{Gom}.
\end{Ex}

Suppose $Y$ is a homology sphere, which is the boundary of a
Stein surface $(X, J)$ with $J$ the associated complex structure on $X$ ,
 then $Y$ has a natural induced holomorphically
fillable contact structure $\xi$ (the $J$-invariant tangent
plane distributions). Define $\mu (Y, \xi) = SW_{X, J}$, this
definition can be extended to the symplectically fillable contact 
structure $\xi$ on $Y$. 

\begin{Pro}
Let $J_1, J_2$ be two complex structures on a Stein surface $X$,
if the induced contact structure $\xi_1, \xi_2$ on $Y$ is isotopic,
that is, there is a   diffeomorphism $f: Y \to Y$ which is homotopic to
identity on $Y$ and sends $\xi_1$ to $\xi_2$, then 
$\mu (Y, \xi_1) = \mu (Y, \xi_2)$.
\end{Pro}
 
\begin{proof}
Under these assumptions, Lisca and Matic prove that $c_1(J_1) = c_1(J_2)$
hence $J_1, J_2$ induce the same $Spin^C$-structure. From the
definition of the relative Seiberg-Witten invariant for $(Y, \xi_i)$, one
 can obtain $\mu (Y, \xi_1) = \mu (Y, \xi_2)$.
\end{proof} 
  
Since any homology 3-sphere $Y$  is obtained by surgery along  
a link $L = \bigcup_{i=1}^n K_i$ in $S^3$, whose   link matrix 
 $Q = (a_{ij})$ (where $a_{ij} $ is the link number for $K_i$ and $K_j$
when $j \ne j$, and $a_{ii}$ is the framing on $K_i$) 
 is  unimodular and  symmetric. One can get a four-manifold $X$ bounded
by $Y$ by attaching 2-handles $h_i = D^2 \times D^2$ ($i=1, \cdots, n$) along 
each knot $K_i$ in $S^3 = \partial ( D^4) $ with framing given by $a_{ii}$,
that is, $ X= D^4 \bigcup_{L} (h_L)$. Note that $H_2( X, {\Bbb Z})$ has 
a basis $\{ \alpha_1, \cdots, \alpha_n \}$, determined  by $\{ K_1, \cdots,
K_n \}$, such that the   intersection matrix is $Q$. 
From the works of Eliashberg and Gompf, we know that if the link $L$
can be realized 
by a Legendrian link in $S^3$ for the canonical contact structure,
 such that  the Thurston-Bennequin invariant for each knot $K_i$ is 
$ a_{ii} + 1 $, then $X$ has a Stein structure $J$. Moreover, the first 
Chern class  $c_1(J) \in H^2( X, {\Bbb Z}) $  of such a Stein structure
$J$ is represented by a cocycle whose value   on  each    such basis element
 $\alpha_i$ is the
rotation number of the corrresponding oriented Legendrian link
component. We can apply  the
relative Seiberg-Witten invariants to study these induced holomorphically
fillable contact structures $\xi_J$ on $Y$, simple examples
show that these relative Seiberg-Witten invariants can distinguish
certain induced contact structures.  Full understanding
of these relative invariants will be interesting  for further
investigations.

\noindent{\bf Acknowledgement:} I wish to acknowledge my 
gratitude to my supervisor Alan Carey for his advice, encouragement
and support throughout this work. I am grateful to Matilde Marcolli,
Tom Mrowka,  Ruibin Zhang for many invaluable discussions.

\vspace{2cm}

\noindent  Department of Pure Mathematics  

\noindent University of Adelaide

\noindent Adelaide, SA 5005

\noindent Australia

\noindent bwang@maths.adelaide.edu.au
\end{document}